\documentclass[11pt]{article}
\usepackage{a4wide,color,comment}

\usepackage{amssymb,amsfonts,amsmath}

\newcommand{\aut}{{\rm Aut}}
\newcommand{\x}{\chi}
\newcommand{\g}{\gamma}

\newcommand{\qed}{\hfill$\Box$}

\renewcommand{\L}{{\rm gr}}

\newcommand{\ia}[1]{{\rm I}_{#1}{\rm A}}

\newcommand{\pf}{{\noindent \it Proof.\ }}

\newtheorem{theorem}{Theorem}
\newtheorem{corollary}{Corollary}
\newtheorem{lemma}{Lemma}
\newtheorem{remark}{Remark}
\newtheorem{proposition}{Proposition}
\newtheorem{example}{Example}

\title{Quotient groups of IA-automorphisms of a free group of rank $3$}
\author{V. Metaftsis, A.I. Papistas and H. Sevaslidou}
\date{}

\begin{document}

\maketitle

\begin{abstract}
We prove that, for any positive integer $c$, the quotient group $\gamma_{c}(M_{3})/\gamma_{c+1}(M_{3})$ of the lower central series of the McCool group $M_{3}$ is isomorphic to two copies of the quotient group $\gamma_{c}(F_{3})/\gamma_{c+1}(F_{3})$ of the lower central series of a free group $F_{3}$ of rank $3$ as $\mathbb{Z}$-modules. Furthermore, we give a necessary and sufficient condition  whether the associated graded Lie algebra ${\rm gr}(M_{3})$ of $M_3$ is naturally embedded into the Johnson Lie algebra ${\cal L}({\rm IA}(F_{3}))$ of the IA-automorphisms of $F_{3}$. 
\end{abstract}


\section{Introduction and Notation}

Let $G$ be a group. For a positive integer $c$, let $\g_c(G)$ be the $c$-th term of
the lower central series of $G$. We point out that $\gamma_{2}(G) = G^{\prime}$; that is, the derived group of $G$. We write ${\rm IA}(G)$ for the kernel of the
natural group homomorphism from ${\rm Aut}(G)$ into ${\rm
Aut}(G/G^{\prime})$ and we call its elements \emph{IA-automorphisms} of $G$. For a
positive integer $c \geq 2$, the natural group epimorphism from
$G$ onto $G/\g_{c}(G)$ induces a group homomorphism $\pi_{c}$ from
 ${\aut}(G)$ of $G$ into 
${\aut}(G/\g_{c}(G))$ of $G/\g_{c}(G)$. Write ${\rm I}_{c}{\rm A}(G) = {\rm
Ker}\pi_{c}$. Note that $\ia{2}(G) = {\rm IA}(G)$. It is proved by
Andreadakis \cite[Theorem 1.2]{andreadakis}, that if $G$ is
residually nilpotent (that is, $\bigcap_{c\geq 1}\gamma_{c}(G) =
\{1\}$), then $\bigcap_{c \geq 2} \ia{c}(G) = \{1\}$. 

Throughout this paper, by \lq \lq Lie algebra\rq \rq,  we mean Lie
algebra over the ring of integers $\mathbb{Z}$. Let $G$ be a
group. Write ${\rm gr}_{c}(G) = \g_c(G)/\g_{c+1}(G)$ for $c \geq 1$ and 
denote by $(a,b)$ the commutator $(a,b)=a^{-1}b^{-1}ab$ with $a, b \in G$. The (restricted) direct sum of
the quotients ${\rm gr}_{c}(G)$ is the {\it associated graded
Lie algebra} of $G$, ${\rm gr}(G) = \bigoplus_{c \geq 1}
{\rm gr}_{c}(G)$. The Lie bracket multiplication in $\L(G)$ is $
[a \gamma_{c+1}(G), b \gamma_{d+1}(G)] = (a, b)
\gamma_{c+d+1}(G)$, 
where $a \gamma_{c+1}(G)$ and $b \gamma_{d+1}(G)$ are the
images of the elements $a \in \gamma_{c}(G)$ and $b \in
\gamma_{d}(G)$ in the quotient groups ${\rm gr}_{c}(G)$ and
${\rm gr}_{d}(G)$, respectively, and $(a, b) \gamma_{c+d+1}(G)$
is the image of the group commutator $(a, b)$ in the quotient
group ${\rm gr}_{c+d}(G)$. Multiplication is then extended to
${\rm gr}(G)$ by linearity.

For a
positive integer $n$, with $n \geq 2$, we write $F_{n}$ for a free
group of rank $n$ with a free generating set $\{x_{1}, \ldots,
x_{n}\}$. For $c \geq 2$, we
write $F_{n,c-1} = F_{n}/\gamma_{c}(F_{n})$. Thus, $F_{n,c-1}$ is the free nilpotent group of rank $n$ and class $c-1$. The natural group epimorphism from $F_{n}$ onto $F_{n,c-1}$
induces a group homomorphism $
\pi_{n,c-1}: {\rm Aut}(F_{n}) \longrightarrow {\rm Aut}(F_{n,c-1})$. 
We write $
{\rm I}_{c}{\rm A}(F_{n}) = {\rm Ker}\pi_{n,c-1}$. It is well known that, for $t, s \geq
2$, $({\rm I}_{t}{\rm A}(F_{n}), {\rm I}_{s}{\rm A}(F_{n})) \subseteq {\rm
I}_{t+s-1}{\rm A}(F_{n})$.  
Since $F_{n}$ is residually nilpotent, we have $\bigcap_{c \geq 2}{\rm I}_{c}{\rm A}(F_{n}) = \{1\}$.  Since $
F_{n,c}/{\rm gr}_{c}(F_{n}) \cong F_{n,c-1}$ 
and ${\rm gr}_{c}(F_{n})$ is a fully invariant subgroup of $F_{n,c}$, the natural group epimorphism from
$F_{n,c}$ onto $F_{n,c-1}$ induces a group homomorphism $
\psi_{c,c-1}: {\rm Aut}(F_{n,c}) \longrightarrow
{\rm Aut}(F_{n,c-1})$. 
It is well known that $\psi_{c,c-1}$ is onto. For $c \geq 2$, we define $
A^{*}_{c}(F_{n}) =
{\rm Im}\pi_{n,c} \cap {\rm Ker}\psi_{c,c-1}$.  
For $t \in \{2, \ldots, c\}$, the natural group epimorphism from $F_{n,c}$ onto $F_{n,c}/\gamma_{t}(F_{n,c})$ induces a group homomorphism $
\sigma_{c,t}: {\rm Aut}(F_{n,c}) \rightarrow {\rm Aut}(F_{n,c}/\gamma_{t}(F_{n,c}))$. Write $
{\rm I}_{t}{\rm A}(F_{n,c}) = {\rm Ker}\sigma_{c,t}$, 
and, for $t = 2$, ${\rm IA}(F_{n,c}) = {\rm I}_{2}{\rm A}(F_{n,c})$. We note that, for $c \geq 2$, $
F_{n,c}/\gamma_{c}(F_{n,c}) \cong F_{n,c-1}$. Thus, for $c \geq 2$,
$A^{*}_{c}(F_{n}) = {\rm Im}\pi_{n,c} \cap {\rm I}_{c}{\rm A}(F_{n,c})$. It is easily proved that $A^{*}_{c}(F_{n}) \cong {\rm I}_{c}{\rm A}(F_{n})/{\rm
I}_{c+1}{\rm A}(F_{n})$ as free abelian groups (see, also, \cite[Section 4, p. 246]{andreadakis}). Furthermore, for $n, c \geq 2$, ${\rm rank}(A^{*}_{c}(F_{n})) \leq n ~{\rm rank}({\rm gr}_{c}(F_{n})) = \frac{n}{c}\sum_{d|c}\mu(d)n^{c/d}$, where $\mu$ is the M\"{o}bius function. We point out that ${\rm rank}({\rm gr}_{c}(F_{n})) = \frac{1}{c}\sum_{d|c}\mu(d)n^{c/d}$ for all $n, c \geq 2$ (see, for example, \cite{mks}). 

From now on, for $r \geq 2$, we write ${\cal L}^{r}({\rm I}{\rm A}(F_{n}))
= {\rm I}_{r}{\rm A}(F_{n})/{\rm I}_{r+1}{\rm
A}(F_{n})$. Form the (restricted) direct sum of the free 
abelian groups ${\cal L}^{r}({\rm I}{\rm A}(F_{n}))$, and
denoted by 
$$
{\cal L}({\rm IA}(F_{n})) = \bigoplus_{r \geq 2}{\cal L}^{r}({\rm I}{\rm A}(F_{n})).
$$ 
It has the structure of a graded Lie algebra with ${\cal L}^{r}({\rm IA}(F_{n}))$ as component of degree $r-1$ in the grading and Lie multiplication given by 
$$
[\phi {\rm I}_{j+1}{\rm A}(F_{n}), \psi {\rm I}_{\kappa+1}{\rm A}(F_{n})] = (\phi^{-1}\psi^{-1}\phi\psi){\rm I}_{j+\kappa}{\rm A}(F_{n}),
$$ 
for all $\phi \in {\rm I}_{j}{\rm A}(F_{n})$, $\psi \in {\rm I}_{\kappa}{\rm A}(F_{n})$ $(j, \kappa \geq 2)$. Multiplication is then extended to
${\cal L}({\rm IA}(F_{n}))$ by linearity. The above Lie algebra is usually called \emph{the
Johnson Lie algebra of} ${\rm IA}(F_{n})$. We point out that, for a positive integer
$c$, $ \gamma_{c}({\rm IA}(F_{n})) \subseteq {\rm I}_{c+1}{\rm
A}(F_{n})$. 

Let $H$ be a finitely generated subgroup of ${\rm IA}(F_{n})$ with $H/H^{\prime}$ torsion-free.  For
a positive integer $q$, let ${\cal L}^{q}_{1}(H) =
\gamma_{q}(H)({\rm I}_{q+2}{\rm A}(F_{n}))/{\rm I}_{q+2}{\rm
A}(F_{n})$. Form the
(restricted) direct sum of abelian groups $
{\cal L}_{1}(H) = \bigoplus_{q \geq 1} {\cal L}^{q}_{1}(H)$.  
It is easily verified that ${\cal
L}_{1}(H)$ is a Lie subalgebra of ${\cal
L}({\rm IA}(F_{n}))$. Furthermore, if $\{y_{1}H^{\prime}, \ldots, y_{m}H^{\prime}\}$ is a $\mathbb{Z}$-basis for $H/H^{\prime}$, then ${\cal L}_{1}(H)$ is generated as Lie algebra by the set $\{y_{1} ({\rm I}_{3}{\rm A}(F_{n})), \ldots, y_{m}({\rm I}_{3}{\rm A}(F_{n}))\}$. By \emph{a natural embedding} of ${\rm gr}(H)$ into ${\cal L}({\rm IA}(F_{n}))$, we mean that there exists a Lie algebra isomorphism $\phi$ from ${\rm gr}(H)$ onto ${\cal L}_{1}(H)$ satisfying the conditions $\phi(y_{i}H^{\prime}) = y_{i}({\rm I}_{3}{\rm A}(F_{n}))$, $i = 1, \ldots, m$. In this case, we also say that ${\rm gr}(H)$ is naturally isomorphic to ${\cal L}_{1}(H)$.  

For $n \geq 2$, it was shown by Magnus \cite{magnus}, using work of
Nielsen \cite{nielsen}, that IA$(F_{n})$ has a finite generating
set $\{\x_{ij}, \x_{ijk}: 1 \leq i, j, k \leq n; i \neq j, k; j <
k\}$, where $\x_{ij}$ maps $x_{i}\mapsto x_{i}(x_{i},x_{j})$ and
$\x_{ijk}$ maps $x_{i}\mapsto x_{i}(x^{-1}_{j},x^{-1}_{k})$, with
both $\x_{ij}$ and $\x_{ijk}$ fixing the remaining basis elements.  Let $M_{n}$ be the subgroup of IA$(F_{n})$ generated by the subset
$S=\{ \x_{ij} : 1\le i,j\le n;\ i\neq j\}$. Then, $M_n$ is called
the {\it McCool group} or the {\it basis conjugating automorphisms
group}. It is easily verified that the following relations are
satisfied by the elements of $S$, provided that, in each case, the
subscripts $i, j, k, q$ occurring are distinct:
$(\chi_{ij}, \chi_{kj}) = 
(\chi_{ij}, \chi_{kq}) =
(\chi_{ij} \chi_{kj}, \chi_{ik}) = 1$. It has been proved in \cite{mccool} that $M_n$ has a presentation
$\langle S\mid Z \rangle$, where $Z$ is the set of all possible
relations of the above forms. Since $\gamma_{c}(M_{n}) \subseteq
\gamma_{c}({\rm IA}(F_{n})) \subseteq {\rm I}_{c+1}{\rm A}(F_{n})$
for all $c \geq 1$, and since $F_{n}$ is residually nilpotent, we
have $\bigcap_{c \geq 1} \gamma_{c}(M_{n}) = \{1\}$ and so,
$M_{n}$ is residually nilpotent. 

In the present paper, we show the following result.

\begin{theorem}\label{th1}
\begin{enumerate}
\item For a positive integer $c$, 
$$
\gamma_{c}(M_{3})/\gamma_{c+1}(M_{3}) \cong \gamma_{c}(F_{3})/ \gamma_{c+1}(F_{3}) \oplus \gamma_{c}(F_{3})/ \gamma_{c+1}(F_{3})
$$
as free abelian groups.  

\item Let $H$ be the subgroup of $M_{3}$ generated by $\chi_{21}, \chi_{12}, \chi_{23}$. Then, ${\cal L}_{1}(M_{3})$ is additively equal to the direct sum of the Lie subalgebras ${\cal L}_{1}(H)$ and ${\cal L}_{1}({\rm Inn}(F_{3}))$, where ${\rm Inn}(F_{3})$ denotes the group of inner automorphisms of $F_{3}$. Furthermore, ${\rm gr}(M_{3})$ is naturally isomorphic to ${\cal L}_{1}(M_{3})$ as Lie algebras if and only if ${\rm gr}(H)$ is naturally isomorphic to ${\cal L}_{1}(H)$ as Lie algebras. 

\end{enumerate}
\end{theorem}     
   
In \cite[Theorem 1]{metpap}, it is shown that $M_{3}$ is a Magnus group. The proof of it was long and tedious. In Section 2, we present a rather simple proof avoiding many of the technical results. The new approach gives us the description of each quotient group $\gamma_{c}(M_{3})/\gamma_{c+1}(M_{3})$ as in Theorem \ref{th1} (1). By a result of Sjogren \cite{sjo} (see, also, \cite[Corollary 1.9]{gupta}), $M_{3}$ satisfies the Subgroup Dimension Problem. That is, each $\gamma_{c}(M_{3})$ is equal to the $c$-th dimension subgroup of $M_{3}$. Furthermore, the new approach helps us to give a necessary and sufficient condition for a natural embedding of  ${\rm gr}(M_{3})$ into ${\cal L}({\rm IA}(F_{3}))$. For $n \geq 2$, let ${\rm Inn}(F_{n})$ denote the subgroup of ${\rm IA}(F_{n})$ consisting of all inner automorphisms of $F_{n}$. In Section 3, by using an observation of Andreadakis \cite[Section 6, p. 249]{andreadakis}, we show that ${\rm gr}({\rm Inn}(F_{n}))$ is naturally embedded into ${\cal L}({\rm IA}(F_{n}))$. Hence, $\gamma_{c}({\rm Inn}(F_{n}))/\gamma_{c+1}({\rm Inn}(F_{n}))$ is isomorphic to a subgroup of ${\cal L}^{c+1}({\rm IA}(F_{n}))$ for all $c \geq 1$. Since ${\rm Inn}(F_{n}) \cong F_{n}$, we obtain 
$$
\frac{1}{c}\sum_{d|c}\mu(d)n^{c/d} \leq {\rm rank}({\cal L}^{c+1}({\rm IA}(F_{n})))
$$
for all $n, c$, with $n \geq 2$. For $c = 1$ and $n \geq 2$, we have ${\rm rank}({\cal L}^{2}({\rm IA}(F_{n}))) = \frac{n^{2}(n-1)}{2}$ (see, \cite[Theorem 5.1]{andreadakis}). For $n = 2$ we have ${\rm IA}(F_{2}) = {\rm Inn}(F_{2})$, by a result of Nielsen \cite{nielsen1} and by a result of Andreadakis \cite[Theorem 6.1]{andreadakis}, we have ${\rm rank}({\cal L}^{c+1}({\rm IA}(F_{2}))) = \frac{1}{c}\sum_{d|c}\mu(d)2^{c/d}$. For $n = 3$, by Theorem \ref{th1}(2), we may give a lower bound of ${\rm rank}({\cal L}^{c+1}({\rm IA}(F_{3})))$ in terms of the rank of ${\cal L}^{c}_{1}(H)$ for all $c$. In fact, we observe that 
$$
\frac{1}{c} \sum_{d|c}\mu(d)2^{c/d} + \frac{1}{c} \sum_{d|c}\mu(d)3^{c/d} \leq {\rm rank}({\cal L}^{c+1}({\rm IA}(F_{3})))
$$  
(see, Remark \ref{re2} below). For $n \geq 4$ and $c \geq 2$, Satoh \cite[Corollary 3.3]{satoh} provides a lower bound for      
${\rm rank}({\cal L}^{c+1}({\rm IA}(F_{n}))).$ 

\section{The associated Lie algebra of $M_{3}$}

\subsection{Lazard elimination}

For a free $\mathbb{Z}$-module $A$, let $L(A)$ be the free
Lie algebra on $A$, that is, the free Lie algebra on $\cal A$, where
$\cal A$ is an arbitrary $\mathbb{Z}$-basis of $A$. Thus, we may
write $L(A) = L({\cal A})$. For a positive integer $c$, let
$L^{c}(A)$ denote the $c$th homogeneous component of $L(A)$. It
is well-known that $
L(A) = \bigoplus_{c \geq 1}L^{c}(A)$. 
For $\mathbb{Z}$-submodules $A$ and $B$ of any Lie algebra, let $[A,B]$ be the $\mathbb{Z}$-submodule spanned by
$[a,b]$ where $a \in A$ and $b \in B$. Furthermore, $B \wr A$
denotes the $\mathbb{Z}$-submodule defined by $
B \wr A = B + [B,A] + [B,A,A] + \cdots$. 

Throughout this paper, we use the left-normed convention for Lie commutators. The following
result is a version of Lazard's "Elimination Theorem" (see
\cite[Chapter 2, Section 2.9, Proposition 10]{bour}). In the form
written here it is a special case of \cite[Lemma 2.2]{bks1} (see, also, \cite[Section 2.2]{metpap}).

\begin{lemma}\label{le1}
Let $U$ and $V$ be free $\mathbb{Z}$-modules, and consider the
free Lie algebra $L(U \oplus V)$. Then, $U$ and $V \wr U$ freely
generate Lie subalgebras $L(U)$ and $L(V \wr U)$, and there is a
$\mathbb{Z}$-module decomposition $L(U \oplus V) = L(U) \oplus L(V
\wr U)$. Furthermore, $
V \wr U = V \oplus [V,U] \oplus [V, U, U] \oplus \cdots$ 
and, for each $n \geq 0$, there is a $\mathbb{Z}$-module isomorphism $
[V, \underbrace{U, \ldots, U}_{n}] \cong V
\otimes \underbrace{U \otimes \cdots \otimes U}_{n}$ under which $[v,u_{1}, \ldots, u_{n}]$ corresponds to $v \otimes u_{1} \otimes \cdots \otimes u_{n}$ for all $v \in V$ and all $u_{1}, \ldots u_{n} \in U$. 
\end{lemma}

As a consequence of Lemma \ref{le1}, we have the following result.

\begin{corollary}\label{c1}
For free $\mathbb{Z}$-modules $U$ and $V$, we write $L(U \oplus V)$ for the free
Lie algebra on $U \oplus V$. Then, there is a
$\mathbb{Z}$-module decomposition $L(U \oplus V) = L(U) \oplus L(V) \oplus L(W)$,
where $W = W_{2} \oplus W_{3} \oplus \cdots$ such that, for all $m
\geq 2$, $W_{m}$ is the direct sum of submodules $[V,U, \underbrace{U, \ldots, U}_{a}, \underbrace{V, \ldots, V}_{b}]$ with $a+b = m-2$ and $a, b \geq 0$. Each $[V,U, \underbrace{U, \ldots, U}_{a}, \underbrace{V, \ldots, V}_{b}]$ is isomorphic to $V \otimes U \otimes \underbrace{U \otimes \cdots \otimes U}_{a} \otimes \underbrace{V \otimes \cdots \otimes V}_{b}$ as
$\mathbb{Z}$-module. Furthermore, $L(W)$ is the ideal of $L(U \oplus V)$ generated by the module
$[V,U]$.
\end{corollary}

\subsection{A decomposition of a free Lie algebra}\label{decomp}

Let $X$ be the free  $\mathbb{Z}$-module with a $\mathbb{Z}$-basis $
\{x_{1}, \ldots, x_{6}\}$ and  $L = L(X)$ the free Lie algebra on $X$. 
For $i = 1, 2, 3$, let $v_{2i} = x_{2i-1}+x_{2i}$. Furthermore, we write 
$$
\begin{array}{ll}
U = {\mathbb{Z}}{\mbox{-}}{\rm span}\{x_{1}, x_{3}, x_{5}\}~~{\rm and}~& V = {\mathbb{Z}}\mbox{-}{\rm span}\{v_{2}, v_{4}, v_{6}\}.
\end{array}
$$
Since $X = U \oplus V$, we have $L = L(U \oplus V)$ and so, $L$ is free on ${\cal X} = \{x_{1}, x_{3}, x_{5}, v_{2}, v_{4}, v_{6}\}$. Let ${\cal J}$ be the subset of $L$,
$$
\begin{array}{ll}
{\cal J}= & \{[v_{2},x_{1}], [v_{4},x_{3}], [v_{6},x_{5}], [v_{4},v_2]-[v_4,x_1], [v_{2},v_4]-[v_2,x_3], [v_{4},v_6]-[v_4,x_5],   \\
& \\
& [v_{6},x_{1}], [v_{6},x_{3}], [v_{2},x_{5}]\}.
\end{array}
$$
The aim of this section is to show the following result. 

\begin{proposition}\label{pr1}
With the above notation, let $L = L(U \oplus V)$ be the free Lie algebra on $U \oplus V$. Let $J$ be the ideal of $L$ generated by the set ${\cal J}$. Then, $L = L(U) \oplus L(V) \oplus J$. Moreover, $J$ is a free Lie algebra.    
\end{proposition}
 
For non negative integers $a$ and $b$, we write $[V,U,~_{a}U, ~_{b}V]$ for $[V,U, \underbrace{U, \ldots, U}_{a}, \underbrace{V, \ldots, V}_{b}]$. By Lemma \ref{le1} and Corollary \ref{c1}, we have 
$$
\begin{array}{lll}
L & = & L(U \oplus V) \\
& = & L(U) \oplus L(V \wr U) \\
& = & L(U) \oplus L(V) \oplus L(W),
\end{array}
$$
where $W = W_{2} \oplus W_{3} \oplus \cdots$ such that, for all $m \geq 2$, 
$$
W_{m} = \bigoplus_{a+b=m-2}[V,U,~_{a}U, ~_{b}V].
$$
Furthermore, $L(V \wr U)$  and $L(W)$ are the ideals in $L$ generated by the modules $V \wr U$ and $[V,U]$, respectively. In particular, $L(W)$ is the ideal in $L$ generated by the natural $\mathbb{Z}$-basis 
$$
[{\cal V}, {\cal U}] = \{[v_{2},x_{1}], [v_{4},x_{3}], [v_{6},x_{5}], [v_{4},x_{1}], [v_{2},x_{3}], [v_{4},x_{5}], [v_{6},x_{1}], [v_{6},x_{3}], [v_{2},x_{5}]\}
$$ 
of $[V,U]$. Let ${\cal X}_{V,U}$ be the natural $\mathbb{Z}$-basis of $V \wr U$. That is, 
$$
{\cal X}_{V,U} = {\cal V} \cup( \bigcup_{a \geq 1}[{\cal V},~_{a}{\cal U}]),
$$
where $[{\cal V},~_{a}{\cal U}]$ is the natural $\mathbb{Z}$-basis of the module $[V,~_{a}U]$. Let $\psi_{2}$ be the $\mathbb{Z}$-linear mapping from $[V,U]$ into $L(V \wr U)$ with
$$
\psi_{2}([v_{4},x_{1}]) = [v_{4},v_{2}] - [v_{4},x_{1}], ~\psi_{2}([v_{2},x_{3}]) = [v_{2},v_{4}] - [v_{2},x_{3}], ~\psi_{2}([v_{4},x_{5}]) = [v_{4},v_{6}] - [v_{4},x_{5}]
$$
and $\psi_{2}$ fixes the remaining elements of $[{\cal V}, {\cal U}]$. It is clear enough that $\psi_{2}$ is a $\mathbb{Z}$-linear monomorphism of $[V,U]$ into $L(V \wr U)$. For $a \geq 3$, let $\psi_{a}$ be the mapping from $[V,U, ~_{(a-2)}U]$ into $L(V \wr U)$ satisfying the conditions $\psi_{a}([v,u,u_{3}, \ldots, u_{a}]) = [\psi_{2}([v,u]),u_{3}, \ldots, u_{a}]$ for all $v \in {\cal V}$ and $u, u_{3}, \ldots, u_{a} \in {\cal U}$. We define a map 
$$
\Psi: {\cal X}_{V,U} \rightarrow L(V \wr U)
$$
by $\Psi(v) = v$ for all $v \in {\cal V}$ and, for $a \geq 2$, $\Psi(v) = \psi_{a}(v)$ for all $v \in [{\cal V}, {\cal U}, ~_{(a-2)}{\cal U}]$. Since $L(V \wr U)$ is free on ${\cal X}_{V,U}$, we obtain $\Psi$ is a Lie algebra homomorphism. By applying Lemma 2.1 in \cite{bks2}, we see that $\Psi$ is a Lie algebra automorphism of $L(V \wr U)$. Since $L(W)$ is a free Lie subalgebra of $L(V \wr U)$ and $\Psi$ is an automorphism, we have $\Psi(L(W))$ is a free Lie subalgebra of $L(V \wr U)$. In fact, 
$$
\Psi(L(W)) = L(\Psi(W)),
$$  
that is, $\Psi(L(W))$ is a free Lie algebra on $\Psi(W)$. 

\begin{lemma}\label{le3}
With the above notation, $L(\Psi(W))$ is an ideal in $L$.
\end{lemma}

\pf Since $\Psi$ is an automorphism of $L(V \wr U)$, we obtain $\Psi(L(W)) = L(\Psi(W))$ is an ideal in $L(V \wr U)$. We point out that 
$$
\begin{array}{lll}
L(V \wr U) & = & \Psi(L(V \wr U)) \\
({\rm By~ Corollary}~ \ref{c1}) & = & \Psi(L(V) \oplus L(W)) \\
(\Psi~{\rm automorphism}) & = & \Psi(L(V)) \oplus \Psi(L(W)) \\
& = & L(V) \oplus L(\Psi(W))
\end{array}
$$
and so, 
$$
\begin{array}{lll}
L & = & L(U) \oplus L(V) \oplus L(\Psi(W)).
\end{array}
$$
To show that $L(\Psi(W))$ is an ideal in $L$, it is enough to show that $[w,u] \in L(\Psi(W))$ for all $w \in L(\Psi(W))$ and $u \in L$. Since $L(\Psi(W))$ is an ideal in $L(V \wr U)$ and Lemma \ref{le1}, it is enough to show that $[w,u] \in L(\Psi(W))$ for all $w \in L(\Psi(W))$ and $u \in L(U)$. Furthermore, we may show that $[w, x_{i_{1}}, \ldots, x_{i_{k}}] \in L(\Psi(W))$ for all $w \in L(\Psi(W))$ and $x_{i_{1}}, \ldots, x_{i_{k}} \in \{x_{1}, x_{3}, x_{5}\}$. Write 
$$
{\cal C} = \Psi([{\cal V}, {\cal U}]) \cup (\bigcup_{a+b \geq 1 \atop a, b \geq 0}[\Psi([{\cal V}, {\cal U}]), ~_{a}{\cal U}, ~_{b}{\cal V}]). \eqno(1)
$$
Since $\cal C$ is a $\mathbb{Z}$-basis for $\Psi(W)$, we have $L(\Psi(W)) = L({\cal C})$. Thus, we need to show that $[w, x_{i_{1}}, \ldots, x_{i_{k}}] \in L(\Psi(W))$ for all $w \in {\cal C}$ and $x_{i_{1}}, \ldots, x_{i_{k}} \in \{x_{1}, x_{3}, x_{5}\}$. Since $\cal C$ is a free generating set of $L(\Psi(W))$, the equation (1) and the linearity of the Lie bracket, we may assume that $w \in [\Psi([{\cal V}, {\cal U}]), ~_{a}{\cal U}, ~_{b}{\cal V}]$ with $a + b \geq 1$. Clearly, we may assume that $b \geq 1$. Since $L(\Psi(W))$ is ideal in $L(V \wr U)$ and, by the Jacobi identity, we may further assume that $w$ has a form $[v, y_{j_{1}}, \ldots, y_{j_{a}}, ~_{\mu}v_{2}, ~_{\nu}v_{4}, ~_{\rho}v_{6}]$ with $y_{j_{1}}, \ldots, y_{j_{a}} \in {\cal U}$, $\mu, \nu, \rho \geq 0$ and $\mu + \nu + \rho \geq 1$. By using the Jacobi identity in the expression $[w, x_{i_{1}}, \ldots, x_{i_{k}}]$, and replacing $[v_{4},x_{1}], [v_{2},x_{3}]$ and $[v_{4},x_{5}]$ by $[v_{4},v_{2}]-\psi_{2}([v_{4},x_{1}])$, $[v_{2},v_{4}]-\psi_{2}([v_{2},x_{3}])$ and $[v_{4},v_{6}]-\psi_{2}([v_{4},x_{5}])$, respectively, as many times as it is needed and since $L(\Psi(W))$ is an ideal in $L(V \wr U)$, we may show that $[w, x_{i_{1}}, \ldots, x_{i_{k}}] \in L(\Psi(W))$. Therefore, $L(\Psi(W))$ is an ideal in $L$. \qed        

\begin{example}\upshape{In the present example, we explain the procedure described in the above proof. Let $w = [\psi_{2}([v_{4},x_{1}]), x_{3}, v_{2}, v_{4}]$. Then, 
$$
\begin{array}{lll}
[w,x_{1}] & = & [\psi_{2}([v_{4},x_{1}]), x_{3}, v_{2}, v_{4}, x_{1}] \\
& = & [\psi_{2}([v_{4},x_{1}]), x_{3}, v_{2}, x_{1}, v_{4}] + [\psi_{2}([v_{4},x_{1}]), x_{3}, v_{2}, [v_{4}, x_{1}]] \\
& = & [\psi_{2}([v_{4},x_{1}]), x_{3}, v_{2}, x_{1}, v_{4}] + \\
& & [\psi_{2}([v_{4},x_{1}]), x_{3}, v_{2}, [v_{4},v_{2}]] - [\psi_{2}([v_{4},x_{1}]), x_{3}, v_{2}, \psi_{2}([v_{4},x_{1}])] \\ 
& = & [\psi_{2}([v_{4},x_{1}]), x_{3}, x_{1}, v_{2}, v_{4}] + [\psi_{2}([v_{4},x_{1}]), x_{3}, [v_{2},x_{1}], v_{4}] +\\
& & [\psi_{2}([v_{4},x_{1}]), x_{3}, v_{2}, [v_{4},v_{2}]] - [\psi_{2}([v_{4},x_{1}]), x_{3}, v_{2}, \psi_{2}([v_{4},x_{1}])] \in L(\Psi(W)).
\end{array}
$$  }
\end{example}

\noindent\emph{Proof of Proposition \ref{pr1}.} Since $\psi_{2}([V,U]) \subseteq J$ and $J$ is ideal in $L$, 
we get $
L(\Psi(W)) \subseteq J$.  
Since ${\cal J} \subseteq L(\Psi(W))$, by Lemma \ref{le3}, we obtain $J \subseteq L(\Psi(W))$. 
Therefore, $J = L(\Psi(W))$. That is, $J$ is a free Lie algebra.   
Furthermore, $
L = L(U) \oplus L(V) \oplus J$. \qed

\vskip .120 in

For $c \geq 2$, let $J^{c} = J \cap L^{c}$. Since $J$ is homogeneous, we get $J = \bigoplus_{c \geq 2}J^{c}$. From the above proof, we have the following result. 

\begin{corollary}\label{co2}
With the above notation, let $\Psi$ be the Lie algebra automorphism of $L(V \wr U)$ defined naturally on $V \wr U$ by means of $\psi_{2}$. Then, $J = L(\Psi(W))$. Furthermore, for $c \geq 2$, $L^{c} = L^{c}(U) \oplus L^{c}(V) \oplus J^{c}$.  
\end{corollary}

\subsection{A description of ${\rm gr}(M_{3})$}

Our aim in this section is to show the following result. For its proof, we use similar arguments as in \cite[Section 6]{metpap}. 

\begin{theorem}\label{th2}
Let $M_{3}$ be the McCool subgroup of the IA-automorphisms ${\rm IA}(F_{3})$ of $F_{3}$. Then, ${\rm gr}(M_{3}) \cong L/J$ as Lie algebras. In particular, ${\rm gr}(M_{3})$ is isomorphic as a Lie algebra to an (external) direct sum of two copies of a free Lie algebra of rank $3$. Furthermore, for each $c$, 
$$
\gamma_{c}(M_{3})/\gamma_{c+1}(M_{3}) \cong \gamma_{c}(F_{3})/\gamma_{c+1}(F_{3}) \oplus \gamma_{c}(F_{3})/\gamma_{c+1}(F_{3})
$$ 
as free abelian groups.
\end{theorem}

Following the notation of the previous subsection, let us denote by $F$ the free group generated by $\{x_1,\ldots, x_6\}$. It is well
known that $\L(F)$ is a free Lie algebra of rank $6$; freely
generated by the set $\{x_{i}F^{\prime}: i = 1, \ldots, 6\}$. The
free Lie algebras $L$ and $\L(F)$ are isomorphic and from now on we identify the two Lie algebras. Furthermore, $L^{c} = \gamma_{c}(F)/\gamma_{c+1}(F)$ for all $c \geq 1$.
Define
$$
\begin{array}{lll}
r_{1} = (x_{1}, x_{2}), & r_{2} = (x_{3}, x_{4}), & r_{3} = (x_{5}, x_{6}) \\
r_{4} = (x_{1}x_{2},x_{5}), & r_{5} = (x_{3}x_{4}, x_{6}), & r_{6} = (x_{1}x_{2}, x_{4}) \\
r_{7} = (x_{3}x_{4}, x_{2}), & r_{8} = (x_{5}x_{6}, x_{3}), &  r_{9}
= (x_{5}x_{6}, x_{1}),
\end{array}
$$
and ${\cal R} = \{r_{1}, \ldots, r_{9}\}$.

Let $N = {\cal R}^{F}$ be the normal closure of $\cal R$ in $F$. For a positive integer $d$, let $N_d=N\cap\g_d(F)$. We point out that for $d \leq 2$, $N_{d} = N$. Further, for $d \geq 2$, $N_{d+1} = N_{d} \cap \gamma_{d+1}(F)$. Define ${\cal I}_{d}(N) = N_{d}\gamma_{d+1}(F)/\gamma_{d+1}(F)$. It is easily verified that ${\cal I}_{d}(N) \cong N_{d}/N_{d+1}$ as $\mathbb{Z}$-modules. The following result was shown in \cite[Section 6]{metpap}.

\begin{proposition}\label{1}
For a positive integer $c$, $N_{c+2}$ is generated by the set
$\{(r^{\pm 1}, g_{1},$ $\ldots, g_{s}): r \in {\cal R}, s \geq c,
g_{1},\ldots, g_{s} \in F \setminus \{1\}\}$. Furthermore, ${\cal
I}_{c+2}(N) = J^{c+2}$ for all $c \geq 1$.
\end{proposition}

Since $F$ is residually nilpotent, we have $ \bigcap_{d \geq
2}N_{d} = \{1\}$. Also $N \subseteq F^{\prime}$, we get ${\cal
I}_{1}(N) = 0$. Moreover, ${\cal I}_{d}(N) \cong N_{d}/N_{d+1}$ as
$\mathbb{Z}$-modules for all $d \geq 2$, and by Proposition
\ref{1} we have, $N_{d} \neq N_{d+1}$ for all $d \geq 2$. Define
$$
{\cal I}(N) = \bigoplus_{d \geq 2} N_{d}
\gamma_{d+1}(F)/\gamma_{d+1}(F) = \bigoplus_{d \geq 2}{\cal
I}_{d}(N).
$$
Since $N$ is a normal subgroup of $F$, we have ${\cal I}(N)$ is an
ideal of $L$ (see \cite{laz}).

\begin{corollary}
${\cal I}(N) = J$.
\end{corollary}

\pf Since $J = \bigoplus_{d \geq 2}J^{d}$ and ${\cal I}_{2}(N) =
J^{2}$, we have from Proposition \ref{1} that ${\cal I}(N) = J$. \qed

\bigskip

\noindent\emph{Proof of Theorem \ref{th2}.} Since $M_{3}/M^{\prime}_{3} \cong F/N F^{\prime} = F/F^{\prime}$,
we have $\L(M_{3})$ is generated as a Lie algebra by the set
$\{\alpha_{i}: i = 1, \ldots, 6\}$ with $\alpha_{i} =
x_{i}M^{\prime}_{3}$. Since $L$ is a free Lie algebra of rank $6$
with a free generating set $\{x_{1}, \ldots, x_{6}\}$, the map
$\zeta$ from $L$ into $\L(M_{3})$ satisfying the conditions
$\zeta(x_{i}) = \alpha_{i}$, $i = 1, \ldots, 6$, extends uniquely
to a Lie algebra homomorphism. Since $\L(M_{3})$ is generated as a
Lie algebra by the set $\{\alpha_{i}: i = 1, \ldots, 6\}$, we have
$\zeta$ is onto. Hence, $L/{\rm Ker}\zeta \cong \L(M_{3})$ as Lie
algebras. By definition, $J \subseteq {\rm Ker} \zeta$, and so
$\zeta$ induces a Lie algebra epimorphism $\overline{\zeta}$ from
$L/J$ onto $\L(M_{3})$. In particular, $\overline{\zeta}(x_{i} + J)
= \alpha_{i}$, $i = 1, \ldots, 6$. Note that $\overline{\zeta}$
induces $\overline{\zeta}_{c}$, say, a $\mathbb{Z}$-linear mapping
from $(L^{c} + J)/J$ onto $\gamma_{c}(M_{3})/\gamma_{c+1}(M_{3})$.
For $c \geq 2$,
$$
\gamma_{c}(M_{3})/\gamma_{c+1}(M_{3}) \cong
\gamma_{c}(F)\gamma_{c+1}(F)N/\gamma_{c+1}(F)N \cong
\gamma_{c}(F)/(\gamma_{c}(F) \cap \gamma_{c+1}(F)N).
$$
Since $\gamma_{c+1}(F) \subseteq \gamma_{c}(F)$, we have by the
modular law,
$$
\gamma_{c}(F)/(\gamma_{c}(F) \cap \gamma_{c+1}(F)N) =
\gamma_{c}(F)/\gamma_{c+1}(F)N_{c}.
$$
But, by Proposition \ref{1}, for $c \geq 3$,
$$
\gamma_{c}(F)/\gamma_{c+1}(F)N_{c} \cong
(\gamma_{c}(F)/\gamma_{c+1}(F))/{\cal I}_{c}(N) \cong L^{c}/J^{c}.
$$
Since ${\cal I}_{2}(N) = J^{2}$, we obtain, for $c \geq 2$, 
$$
\gamma_{c}(F)/\gamma_{c+1}(F)N_{c} \cong
(\gamma_{c}(F)/\gamma_{c+1}(F))/{\cal I}_{c}(N) \cong L^{c}/J^{c}.
$$
Therefore, for $c \geq 1$, 
$$
\gamma_{c}(M_{3})/\gamma_{c+1}(M_{3}) \cong L^{c}/J^{c} \cong
(L^{c})^{*},
$$
by Corollary \ref{co2}, where $(L^{c})^{*} = L^{c}(U) \oplus L^{c}(V)$. Since both $L(U)$ and $L(V)$ are free Lie algebras of rank $3$, we have $L(U) \cong L(V) \cong {\rm gr}(F_{3})$ in a natural way and so, for $c \geq 1$, $L^{c}(U) \cong L^{c}(V) \cong \gamma_{c}(F_{3})/\gamma_{c+1}(F_{3})$ as free abelian groups. Hence, for $c \geq 1$, 
$$
\gamma_{c}(M_{3})/\gamma_{c+1}(M_{3}) \cong \gamma_{c}(F_{3})/\gamma_{c+1}(F_{3}) \oplus \gamma_{c}(F_{3})/\gamma_{c+1}(F_{3})
$$
as free abelian groups and so, ${\rm rank}(\gamma_{c}(M_{3})/\gamma_{c+1}(M_{3})) = {\rm
rank}(L^{c})^{*}$. Since $J = \bigoplus_{c \geq 2}J^{c}$, we have
$(L^{c} + J)/J \cong L^{c}/(L^{c} \cap J) = L^{c}/J^{c} \cong
(L^{c})^{*}$ and so, we obtain ${\rm
Ker}\overline{\zeta}_{c}$ is torsion-free. Since ${\rm
rank}(\gamma_{c}(M_{3})/\gamma_{c+1}(M_{3})) = {\rm
rank}(L^{c})^{*}$, we have ${\rm Ker}\overline{\zeta}_{c} = \{1\}$
and so, $\overline{\zeta}_{c}$ is isomorphism. Since
$\overline{\zeta}$ is epimorphism and each $\overline{\zeta}_{c}$ is
isomorphism, we have $\overline{\zeta}$ is isomorphism. Hence, $L/J
\cong \L(M_{3})$ as Lie algebras. \qed

\section{Embeddings}

In this section, we shall give a criterion for the natural embedding of ${\rm gr}(M_3)$ into ${\cal L}_1({\rm IA}(F_3))$. We shall prove the following result.

\begin{lemma}\label{le3a}
Let $H$ be a finitely generated subgroup of ${\rm IA}(F_{n})$, $n \geq 2$, with $H/H^{\prime}$ torsion-free, and let $\{y_{1}H^{\prime}, \ldots, y_{m}H^{\prime}\}$ be a $\mathbb{Z}$-basis for $H/H^{\prime}$. Then, ${\rm gr}(H)$ is naturally isomorphic to ${\cal L}_{1}(H)$ if and only if $\gamma_{c}(H) = H \cap ({\rm I}_{c+1}{\rm A}(F_{n}))$ for all $c$.
\end{lemma}

\pf We assume that $\gamma_{c}(H) = H \cap ({\rm I}_{c+1}{\rm A}(F_{n}))$ for all $c$. For $c \geq 1$, let $\psi_{c}$ be the natural $\mathbb{Z}$-module epimorphism from ${\rm gr}_{c}(H)$ onto ${\cal L}_{1}^{c}(H)$. Since $\gamma_{c}(H) \cap ({\rm I}_{c+2}{\rm A}(F_{n})) = H \cap ({\rm I}_{c+2}{\rm A}(F_{n})) = \gamma_{c+1}(H)$ for all $c$, we get $\psi_{c}$ is isomorphism for all $c \geq 1$. Since ${\rm gr}(H)$ is the (restricted) direct sum of the quotients ${\rm gr}_{c}(H)$, there exists a group homomorphism $\psi$ from ${\rm gr}(H)$ into ${\cal L}_{1}({\rm IA}(F_{n}))$ such that each $\psi_{c}$ is the restriction of $\psi$ on ${\rm gr}_{c}(H)$. It is easily shown that $\psi$ is a Lie algebra homomorphism. Since $\psi(y_{i}H^{\prime}) = y_{i}({\rm I}_{3}{\rm A}(F_{n}))$, $i = 1, \ldots, n$, we get $\psi$ is a Lie algebra epimorphism. Furthermore, since each $\psi_{c}$ is a $\mathbb{Z}$-module isomorphism, we obtain $\psi$ is a Lie algebra isomorphism. Conversely, let $\phi$ be a Lie algebra isomorphism from ${\rm gr}(H)$ onto ${\cal L}_{1}(H)$ satisfying the conditions $\phi(y_{i}H^{\prime}) = y_{i} ({\rm I}_{3}{\rm A}(F_{n}))$, $i = 1, \ldots, m$. Then, $\phi$ induces a $\mathbb{Z}$-module isomorphism $\phi_{c}$ from ${\rm gr}_{c}(H)$ onto ${\cal L}^{c}_{1}(H)$ for all $c$. In particular, $\phi_{c}((y_{i_{1}}, \ldots, y_{i_{c}})\gamma_{c+1}(H)) = (y_{i_{1}}, \ldots, y_{i_{c}})({\rm I}_{c+2}{\rm A}(F_{n}))$ for all $i_{1}, \ldots, i_{c} \in \{1, \ldots, m\}$. Furthermore, ${\rm gr}_{c}(H)$ is $\mathbb{Z}$-module isomorphic to $\gamma_{c}(H)/(\gamma_{c}(H) \cap ({\rm I}_{c+2}{\rm A}(F_{n})))$. Since ${\rm gr}_{c}(H)$ is polycyclic and so, it is a hopfian group, we have $\gamma_{c+1}(H) = \gamma_{c}(H) \cap ({\rm I}_{c+2}{\rm A}(F_{n}))$ for all $c$. We claim that $\gamma_{c}(H) = H \cap ({\rm I}_{c+1}{\rm A}(F_{n}))$ for all $c$. Since $\gamma_{c}(H) \subseteq {\rm I}_{c+1}{\rm A}(F_{n})$, it is enough to show that $H \cap ({\rm I}_{c+1}{\rm A}(F_{n})) \subseteq \gamma_{c}(H)$. To get a contradiction, let $\alpha \in H \cap ({\rm I}_{c+1}{\rm A}(F_{n}))$ and $\alpha \notin \gamma_{c}(H)$. Since $\gamma_{c}(H) \subseteq {\rm I}_{c+1}{\rm A}(F_{n})$ for all $c$ and $\bigcap_{t \geq 2}{\rm I}_{t}{\rm A}(F_{n}) = \{1\}$, we get $H$ is residually nilpotent. Thus, there exists a unique $d \in {\mathbb{N}}$ such that $\alpha \in \gamma_{d}(H) \setminus \gamma_{d+1}(H)$. Therefore, $\alpha \notin H \cap ({\rm I}_{d+2}{\rm A}(F_{n}))$. Since $\alpha \in \gamma_{d}(H) \setminus \gamma_{d+1}(H)$ and $\alpha \notin \gamma_{c}(H)$, we have $\gamma_{c}(H) \subseteq \gamma_{d+1}(H)$ and so, $d+1 \leq c$. Let $k$ be a non-negative integer such that $c = d+1+k$. Since $H \cap ({\rm I}_{c+1}{\rm A}(F_{n})) = H \cap ({\rm I}_{d+2+k}{\rm A}(F_{n}))$, we have $\alpha \in {\rm I}_{d+2+k}{\rm A}(F_{n})$. But ${\rm I}_{d+2+k}{\rm A}(F_{n}) \subseteq {\rm I}_{d+2}{\rm A}(F_{n})$ and so, $\alpha \in {\rm I}_{d+2}{\rm A}(F_{n})$, which is a contradiction. Therefore, $\gamma_{c}(H) = H \cap ({\rm I}_{c+1}{\rm A}(F_{n}))$ for all $c$.   \qed

\begin{remark}\upshape{It is known that ${\rm IA}(F_{n})/\gamma_{2}({\rm IA}(F_{n}))$, with $n \geq 2$, is torsion-free and its rank is $\frac{n^{2}(n-1)}{2}$. Furthermore, $\gamma_{2}({\rm IA}(F_{n})) = {\rm I}_{3}{\rm A}(F_{n})$ (see, for example, \cite{pet}). It was conjectured by Andreadakis \cite{andreadakis} that $\gamma_{c}({\rm IA}(F_{n})) = {\rm I}_{c+1}{\rm A}(F_{n})$ for all $c$. By Lemma \ref{le3a}, ${\rm gr}({\rm IA}(F_{n}))$ is naturally isomorphic to ${\cal L}_{1}({\rm IA}(F_{n}))$ if and only if Andreadakis conjecture is valid. Now, if $H$ is a finitely generated subgroup of ${\rm IA}(F_n)$ with $H/H'$ torsion free, then the statement ${\rm gr}(H)$ is naturally isomorphic to ${\cal L}_1(H)$ seems to be an ``Andreadakis Conjecture'' for $H$.
}
\end{remark}
 
\subsection{The associated Lie algebra of ${\rm Inn}(F_{n})$}

In this section, we show that the associated Lie algebra ${\rm gr}({\rm Inn}(F_{n}))$ of the inner automorphisms ${\rm Inn}(F_{n})$ of $F_{n}$, with $n \geq 2$, is naturally embedded into ${\cal L}({\rm IA}(F_{n}))$. Throughout this section, we write $E_{n} = {\rm Inn}(F_{n})$. Recall that, for $g \in F_{n}$, $\tau_{g}(x) = gxg^{-1}$ for all $x \in F_{n}$. Thus, $E_{n} = \{\tau_{g}: g \in F_{n}\}$. The following result has been proved in \cite[Section 6]{andreadakis}.

\begin{lemma}\label{le4} For a positive integer $c$, $\gamma_{c}(E_{n}) = E_{n} \cap {\rm I}_{c+1}{\rm A}(F_{n})$.
\end{lemma}

Using the above we may show the following.

\begin{proposition}\label{pr3}
Let $n$ be positive integer, with $n \geq 2$. Then, ${\rm gr}(E_{n})$ 
is naturally embedded into ${\cal L}({\rm IA}(F_{n}))$. In particular, for all $c$, ${\rm gr}_{c}(E_{n})$ is isomorphic to a $\mathbb{Z}$-submodule of ${\cal L}^{c+1}({\rm IA}(F_{n}))$.
\end{proposition}

\pf Since $F_{n}$ is centerless, we have $F_{n} \cong E_{n}$ in a natural way and so, $E_{n}$ is finitely generated. Moreover, ${\rm gr}(E_{n})$ is a free Lie algebra of rank $n$. Namely, ${\rm gr}(E_{n})$ is freely generated by the set $\{\tau_{x_{i}}E^{\prime}_{n}: i = 1, \ldots, n\}$. Let $\phi$ be the mapping from  $\{\tau_{x_{i}}E^{\prime}_{n}: i = 1, \ldots, n\}$ to ${\cal L}_{1}(E_{n})$ satisfying the conditions $\phi(\tau_{x_{i}}E^{\prime}_{n}) = \tau_{x_{i}}{\rm I}_{3}{\rm A}(F_{n})$, $i = 1, \ldots, n$. Since ${\rm gr}(E_{n})$ is free on $\{\tau_{x_{i}}E^{\prime}_{n}: i = 1, \ldots, n\}$, $\phi$ is extended to a Lie algebra epimorphism. By Lemma \ref{le3a}, it is enough to show that $\gamma_{c}(E_{n}) = E_{n} \cap {\rm I}_{c+1}{\rm A}(F_{n})$ for all $c$, which is valid by Lemma \ref{le4}.  
Therefore, ${\rm gr}(E_{n})$ 
is naturally embedded into ${\cal L}({\rm IA}(F_{n}))$. Since 
$$
\begin{array}{lll}
{\cal L}^{c}_{1}(E_{n}) & \cong & \gamma_{c}(E_{n})/(\gamma_{c}(E_{n}) \cap ({\rm I}_{c+2}{\rm A}(F_{n}))) \\ 
& = & \gamma_{c}(E_{n})/(E_{n} \cap {\rm I}_{c+2}{\rm A}(F_{n}))) \\
& = & {\rm gr}_{c}(E_{n})
\end{array}
$$ 
for all $c$, we have ${\rm gr}_{c}(E_{n})$ is isomorphic to a $\mathbb{Z}$-submodule of ${\cal L}^{c+1}({\rm IA}(F_{n}))$ for all $c$. \qed 

\begin{corollary}\label{}
For positive integers $n$ and $c$, with $n \geq 2$, 
$$
\frac{1}{c} \sum_{d|c} \mu(d) n^{c/d} \leq {\rm rank}({\cal L}^{c+1}({\rm IA}(F_{n})))
$$ where $\mu$ is the M\"{o}bius function.
\end{corollary}

\pf Since $F_{n} \cong E_{n}$ as groups and ${\rm gr}(F_{n})$ is a free Lie algebra, we get ${\rm gr}(F_{n}) \cong {\rm gr}(E_{n})$ as Lie algebras. Hence, for $c \geq 1$, ${\rm gr}_{c}(F_{n}) \cong {\rm gr}_{c}(E_{n})$ as $\mathbb{Z}$-modules. Since ${\rm rank}({\rm gr}_{c}(F_{n})) = \frac{1}{c} \sum_{d|c} \mu(d) n^{c/d}$,  we obtain, by Proposition \ref{pr3}, the required result.   \qed

\subsection{The Lie algebra ${\cal L}_{1}(M_{3})$}

Write $b_{1} = \chi_{21}$, $b_{2} = \chi_{12}$, $b_{3} = \chi_{23}$, $u_{2} = \chi_{31}\chi_{21}$, $u_{4} = \chi_{32}\chi_{12}$ and $u_{6} = \chi_{23}\chi_{13}$. Then, $M_{3}$ is generated by the set $\{b_{1}, b_{2}, b_{3}, u_{2}, u_{4}, u_{6}\}$ and its presentation is given by 
\begin{flushleft}
$M_{3} = \langle b_{1}, b_{2}, b_{3}, u_{2}, u_{4}, u_{6}: (u_{2},b_{1}), (u_{2},u_{4}b^{-1}_{2}), (u_{2},b_{3}),$ 
\end{flushleft} 
\begin{flushright}
$(u_{4},u_{2}b^{-1}_{1}), (u_{4},b_{2}), (u_{4},b^{-1}_{3}u_{6}), (u_{6},b_{1}), (u_{6},b_{2}), (u_{6},b_{3}) \rangle.$
\end{flushright}
We point out that $u_{2} = \tau_{x^{-1}_{1}}$, $u_{4} = \tau_{x^{-1}_{2}}$ and $u_{6} = \tau_{x^{-1}_{3}}$. Let $H$ and $E$ be the subgroups of $M_{3}$ generated by the sets $\{b_{1}, b_{2}, b_{3}\}$ and $\{u_{2}, u_{4}, u_{6}\}$, respectively. We point out that $E = {\rm Inn}(F_{3})$. One can easily see that $H$ is a free group of rank $3$. Thus, ${\rm gr}(H) \cong {\rm gr}(E) \cong {\rm gr}(F_{3})$ as Lie algebras.  

\begin{proposition}\label{pr4}
${\cal L}_{1}(M_{3})$ is additively equal to the direct sum of the Lie subalgebras ${\cal L}_{1}(H)$ and ${\cal L}_{1}(E)$.   
\end{proposition}

\pf By the proof of Proposition \ref{pr3}, ${\rm gr}(E) \cong {\cal L}_{1}(E)$ and so, ${\cal L}_{1}(E)$ is a free Lie algebra of rank $3$. Since $b_{1} \notin {\rm I}_{3}{\rm A}(F_{3})$, we have ${\cal L}_{1}(H)$ is a non-trivial subalgebra of ${\cal L}({\rm IA}(F_{3}))$. In fact, ${\cal L}_{1}(H)$ is generated by the set $\{b_{i}({\rm I}_{3}{\rm A}(F_{3})): i = 1, 2, 3\}$. Since $(\tau_{g}, \phi) = \tau_{g^{-1}\phi^{-1}(g)}$ for all $\phi \in {\rm Aut}(F_{3})$ and $g \in F_{3}$, we have ${\cal L}_{1}(E)$ is an ideal in ${\cal L}({\rm IA}(F_{3}))$ and so, ${\cal L}_{1}(H) + {\cal L}_{1}(E)$ is a Lie subalgebra of ${\cal L}_{1}(M_{3})$. Let $\overline{w} \in {\cal L}_{1}(H) \cap {\cal L}_{1}(E)$. Since both ${\cal L}_{1}(H)$ and ${\cal L}_{1}(E)$ are graded Lie algebras, we may assume that $\overline{w} \in {\cal L}^{d}_{1}(H) \cap {\cal L}^{d}_{1}(E)$ for some $d$. Thus, there are $u \in \gamma_{d}(H)$ and $v \in \gamma_{d}(E)$ such that $\overline{w} = u ({\rm I}_{d+2}{\rm A}(F_{3})) = v ({\rm I}_{d+2}{\rm A}(F_{3}))$. To get a contradiction, we assume that $u, v \notin {\rm I}_{d+2}{\rm A}(F_{3})$. Therefore, $v \in \gamma_{d}(E) \setminus \gamma_{d+1}(E)$ and so, there exists $\omega \in \gamma_{d}(F_{3}) \setminus \gamma_{d+1}(F_{3})$ such that $v = \tau_{\omega} \rho$, where $\rho \in \gamma_{d+1}(E)$. Since $\gamma_{d+1}(E) \subseteq {\rm I}_{d+2}{\rm A}(F_{3})$, we get $u^{-1}\tau_{\omega} \in {\rm I}_{d+2}{\rm A}(F_{3})$. Since $u^{-1}$ fixes $x_{3}$, we have $x^{-1}_{3} (u^{-1}\tau_{\omega}(x_{3})) = (x_{3}, u^{-1}(\omega^{-1})) \in \gamma_{d+2}(F_{3})$. Since $u^{-1}(x_{j}) = x_{j}y_{j}$, $y_{j} \in \gamma_{d+1}(F_{3})$, $j = 1, 2$, and $\gamma_{d}(F_{3})$ is a fully invariant subgroup of $F_{3}$, we have $u^{-1}(\omega^{-1}) = \omega^{-1} \omega_{1}$, with $\omega_{1} \in \gamma_{d+1}(F_{3})$ and so, $(w, x_{3}) \in \gamma_{d+2}(F_{3})$. Since ${\rm gr}(F_{3})$ is a free Lie algebra of rank $3$ with a free generating set $\{x_{i}F^{\prime}_{3}: i = 1,2, 3\}$ and $\gamma_{d}(F_{3})/\gamma_{d+1}(F_{3})$ is the $d$-th homogeneous component of ${\rm gr}(F_{3})$, we obtain $(w, x_{3}) \in \gamma_{d+1}(F_{3}) \setminus \gamma_{d+2}(F_{3})$, which is a contradiction. Therefore, ${\cal L}_{1}(H) \cap {\cal L}_{1}(E) = \{0\}$. By the proof of Theorem \ref{th2}, we obtain ${\cal L}_{1}(M_{3}) = {\cal L}_{1}(H) \oplus {\cal L}_{1}(E)$. \qed

\begin{remark}\label{re2}\upshape{By Proposition \ref{pr4}, for all $c$, we obtain ${\cal L}^{c}_{1}(M_{3}) = {\cal L}^{c}_{1}(H) \oplus {\cal L}^{c}_{1}(E)$. Since ${\cal L}^{c}_{1}(E) \cong \gamma_{c}(E)/\gamma_{c+1}(E) \cong \gamma_{c}(F_{3})/\gamma_{c+1}(F_{3})$, we have ${\rm rank}({\cal L}^{c}_{1}(E)) = \frac{1}{c} \sum_{d|c}\mu(d)3^{c/d}$. Thus, for any $c$, 
$$
{\rm rank}({\cal L}^{c}_{1}(H)) + \frac{1}{c} \sum_{d|c}\mu(d)3^{c/d} \leq {\rm rank}({\cal L}^{c+1}({\rm IA}(F_{3}))).
$$
Let $H_{1}$ be the subgroup of $H$ generated by the set $\{\chi_{21}, \chi_{23}\}$. We point out that $H_{1}$ is a free group of rank $2$. Then, ${\rm gr}(H_{1})$ is a free Lie algebra of rank $2$. Since both $\chi_{21}$ and $\chi_{23}$ fix $x_{1}$ and $x_{3}$, it may be shown that $\gamma_{c}(H_{1}) = H_{1} \cap {\rm I}_{c+1}{\rm A}(F_{3})$ for all $c$. Therefore, $\gamma_{c}(H_{1}) \cap {\rm I}_{c+2}{\rm A}(F_{3}) = \gamma_{c+1}(H_{1})$ for all $c$ and so, ${\cal L}^{c}_{1}(H_{1}) \cong {\rm gr}_{c}(H_{1})$ for all $c \geq 1$. Since ${\rm gr}(H_{1})$ is a free Lie algebra on ${\cal H}_{1} = \{\chi_{21}H^{\prime}_{1}, \chi_{23}H^{\prime}_{1}\}$, the mapping $\psi$ from ${\cal H}_{1}$ into ${\cal L}_{1}(H_{1})$ satisfying the conditions $\psi(aH^{\prime}_{1}) = a({\rm I}_{3}{\rm A}(F_{3}))$ with $a \in \{\chi_{21}, \chi_{23}\}$ can be extended to a Lie algebra homomorphism $\overline{\psi}$. Since ${\cal L}_{1}(H_{1})$ is generated as Lie algebra by the set $\{\chi_{21}({\rm I}_{3}{\rm A}(F_{3})), \chi_{23}({\rm I}_{3}{\rm A}(F_{3}))\}$, we have $\overline{\psi}$ is onto. By Lemma \ref{le3a}, $\overline{\psi}$ is a natural Lie algebra isomorphism from ${\rm gr}(H_{1})$ onto ${\cal L}_{1}(H_{1})$. Since ${\cal L}^{c}_{1}(H_{1}) \cong {\rm gr}_{c}(H_{1})$ for all $c \geq 1$, we have ${\rm rank}({\cal L}^{c}_{1}(H_{1})) = \frac{1}{c} \sum_{d|c}\mu(d)2^{c/d}$. Since ${\cal L}_{1}(H_{1})$ is a Lie subalgebra of ${\cal L}_{1}(H)$, we obtain ${\cal L}^{c}_{1}(H_{1}) \leq {\cal L}^{c}_{1}(H)$ for all $c$. Therefore,           
$$ 
\frac{1}{c} \sum_{d|c}\mu(d)2^{c/d} \leq {\rm rank}({\cal L}^{c}_{1}(H))
$$
for all $c$ and so,
$$
\frac{1}{c} \sum_{d|c}\mu(d)2^{c/d} + \frac{1}{c} \sum_{d|c}\mu(d)3^{c/d} \leq {\rm rank}({\cal L}^{c+1}({\rm IA}(F_{3})))
$$  
for all $c$.
}
\end{remark}

In our next result, we give a necessary and sufficient condition for a natural embedding of ${\rm gr}(M_{3})$ into ${\cal L}({\rm IA}(F_{3}))$. 

\begin{proposition}\label{pr5}
Let $H$ be the subgroup of $M_{3}$ generated by $\chi_{21}, \chi_{12}, \chi_{23}$. Then, ${\rm gr}(M_{3})$ is naturally isomorphic to ${\cal L}_{1}(M_{3})$ as Lie algebras if and only if ${\rm gr}(H)$ is naturally isomorphic to ${\cal L}_{1}(H)$ as Lie algebras.  
\end{proposition}

\pf Suppose that ${\rm gr}(M_{3})$ is naturally isomorphic to ${\cal L}_{1}(M_{3})$ as Lie algebras. Thus, for all $c$, ${\rm gr}_{c}(M_{3}) \cong {\cal L}^{c}_{1}(M_{3})$ as $\mathbb{Z}$-modules. By Proposition \ref{pr4}, ${\rm gr}_{c}(M_{3}) \cong {\cal L}^{c}_{1}(H) \oplus {\cal L}^{c}_{1}(E)$. Since, for all $c$, ${\rm gr}_{c}(M_{3}) \cong {\rm gr}_{c}(F_{3}) \oplus {\rm gr}_{c}(F_{3})$ as $\mathbb{Z}$-modules, we obtain
$$
\begin{array}{lll}
{\rm rank}({\cal L}^{c}_{1}(H)) & = & {\rm rank}({\cal L}^{c}_{1}(E)) \\
 & = & {\rm rank}({\rm gr}_{c}(F_{3})) \\
 & = & {\rm rank}({\rm gr}_{c}(H)).
 \end{array}
 $$
Therefore, ${\cal L}^{c}_{1}(H) \cong {\rm gr}_{c}(H)$ for all $c$. Hence, $\gamma_{c}(H) \cap {\rm I}_{c+2}{\rm A}(F_{3}) = \gamma_{c+1}(H)$ for all $c$ and so, as in the proof of Lemma \ref{le3a}, $\gamma_{c}(H) = H \cap {\rm I}_{c+1}{\rm A}(F_{3})$ for all $c$. We point out that $H/H^{\prime}$ is torsion-free with a free generating set ${\cal H} = \{\chi_{21}H^{\prime}, \chi_{12}H^{\prime}, \chi_{23}H^{\prime}\}$. Since ${\rm gr}(H)$ is a free Lie algebra with a free generating set $\cal H$, and ${\cal L}_{1}(H)$ is generated as a Lie algebra by the set $\{\chi_{21}({\rm I}_{3}{\rm A}(F_{3})), \chi_{12}({\rm I}_{3}{\rm A}(F_{3})), \chi_{23}({\rm I}_{3}{\rm A}(F_{3}))\}$, we have there exists a natural Lie epimorphism from ${\rm gr}(H)$ onto ${\cal L}_{1}(H)$. By Lemma \ref{le3a}, ${\rm gr}(H)$ is naturally isomorphic to ${\cal L}_{1}(H)$ as Lie algebras. 

Conversely, let ${\rm gr}(H)$ be naturally isomorphic to ${\cal L}_{1}(H)$ as Lie algebras. By Lemma \ref{le3a}, $\gamma_{c}(H) = H \cap ({\rm I}_{c+1}{\rm A}(F_{3}))$ for all $c$. Thus, 
$$
\begin{array}{lll}
\gamma_{c}(H) \cap {\rm I}_{c+2}{\rm A}(F_{3}) & = & H \cap {\rm I}_{c+1}{\rm A}(F_{3}) \cap {\rm I}_{c+2}{\rm A}(F_{3}) \\
& = & H \cap {\rm I}_{c+2}{\rm A}(F_{3}) \\
& = & \gamma_{c+1}(H)
\end{array}
$$
and so, ${\cal L}^{c}_{1}(H) \cong {\rm gr}_{c}(H)$ for all $c$. By Proposition \ref{pr4}, 
$$
{\cal L}^{c}_{1}(M_{3}) = {\cal L}^{c}_{1}(H) \oplus {\cal L}^{c}_{1}(E)
$$
for all $c$, and by Theorem \ref{th2}, 
$$
{\cal L}^{c}_{1}(M_{3}) \cong {\rm gr}_{c}(M_{3})
$$
for all $c$. Therefore, 
$$
\gamma_{c}(M_{3}) \cap {\rm I}_{c+2}{\rm A}(F_{3}) = \gamma_{c+1}(M_{3})
$$
for all $c$. As in the proof of Lemma \ref{le3a}, for all $c$, 
$$
M_{3} \cap {\rm I}_{c+1}{\rm A}(F_{3}) = \gamma_{c}(M_{3}).
$$
Let $\psi$ be the Lie algebra epimorphism from $L$ onto ${\cal L}_{1}(M_{3})$ satisfying the conditions $\psi(x_{2j-1}) = b_{j}({\rm I}_{3}{\rm A}(F_{3}))$ and $\psi(v_{2j}) = \tau_{x^{-1}_{2j-1}}({\rm I}_{3}{\rm A}(F_{3}))$, $j = 1, 2, 3$.
By $\tau_w$ we denote the inner automorphisms of $F_3$ with $\tau_w(x)=w^{-1}xw$. Since $J \subseteq {\rm Ker}\psi$ and ${\rm gr}(M_{3}) \cong L/J$ as Lie algebras, there exists a Lie algebra epimorphism $\overline{\psi}$ from ${\rm gr}(M_{3})$ onto ${\cal L}_{1}(M_{3})$ satisfying the conditions $\overline{\psi}(b_{i}M^{\prime}_{3}) = b_{i}({\rm I}_{3}{\rm A}(F_{3}))$ and 
$\overline{\psi}(\tau_{x^{-1}_{i}}M^{\prime}_{3}) = \tau_{x^{-1}_{i}}({\rm I}_{3}{\rm A}(F_{3}))$, $i = 1, 2, 3$. Since, in addition, $M_{3} \cap {\rm I}_{c+1}{\rm A}(F_{3}) = \gamma_{c}(M_{3})$ for all $c$, we obtain from Lemma \ref{le3a} the required result.  \qed

\begin{remark}\label{re3}\upshape{
For $c = 1, 2, 3$, we show that $\gamma_{c}(H) = H \cap {\rm I}_{c+1}{\rm A}(F_{3})$. That is, ${\cal L}^{c}_{1}(H) \cong \gamma_{c}(H)/\gamma_{c+1}(H)$ for $c = 1, 2, 3$. For $c \geq 4$, the method used becomes very cumbersome. For simplicity, let $x = \chi_{21}$, $y = \chi_{12}$ and $z = \chi_{23}$. For $c = 1$, our claim is trivially true. Let $c = 2$. Since $\gamma_{2}(H)/\gamma_{3}(H)$ is a free abelian group of rank $3$, we have each element $h \in \gamma_{2}(H) \setminus \gamma_{3}(H)$ is uniquely written as 
$$
h = (z,x)^{a_{1}}(z,y)^{a_{2}}(y,x)^{a_{3}} v,
$$
where $v \in \gamma_{3}(H)$ and $a_{1}, a_{2}, a_{3}$ are non-negative integers. We claim that $h \notin {\rm I}_{4}{\rm A}(F_{3})$. By a direct calculation, $(z,x)^{a_{1}}(x_{i}) = x_{i}$, $i \neq 2$, and $(z,x)^{a_{1}}(x_{2}) = x_{2}(x_{3},x_{1},x_{2})^{-a_{1}}u$, where $u \in \gamma_{4}(F_{3})$. Working modulo $\gamma_{3}(H)$, we have 
$$
\begin{array}{ccc}
(z,\tau_{x^{-1}_{2}})^{a_{2}} \equiv (z,y)^{a_{2}}(z,\chi_{32})^{a_{2}} & {\rm and} & 
(\tau_{x^{-1}_{2}},x)^{a_{3}} \equiv (y,x)^{a_{3}} (\chi_{32},x)^{a_{3}}. 
\end{array}
$$
We point out that, for $\phi \in {\rm IA}(F_{3})$ and $g \in F_{3}$, we have $(\tau_{g}, \phi) = \tau_{g^{-1}\phi^{-1}(g)}$. Since $(z,\tau_{x^{-1}_{2}}) = \tau_{(x^{-1}_{3},x^{-1}_{2})}$ and $(\tau_{x^{-1}_{2}}, x) = \tau_{(x^{-1}_{2},x^{-1}_{1})}$, we get    
$$
(\tau_{(x^{-1}_{3},x^{-1}_{2})})^{a_{2}} \equiv (z,y)^{a_{2}}(z,\chi_{32})^{a_{2}} \ \ 
{\rm and} \ \  
(\tau_{(x^{-1}_{2},x^{-1}_{1})})^{a_{3}} \equiv (y,x)^{a_{3}} (\chi_{32},x)^{a_{3}}.  \eqno{(*)}
$$
Since both $(z,\chi_{32})^{a_{2}}$ and $(\chi_{32},x)^{a_{3}}$ fix $x_{1}$, we obtain 
$$
\begin{array}{ccc}
(z,y)^{a_{2}}(x_{1}) = x_{1}(x_{1}, (x_{3},x_{2}))^{a_{2}}v_{1} &
{\rm and} & 
(y,x)^{a_{3}}(x_{1}) = x_{1}(x_{1}, (x_{2},x_{1}))^{a_{3}}v_{2},
\end{array}
$$
where $v_{1}, v_{2} \in \gamma_{4}(F_{3})$. Thus, $(z,y), (y,x) \in {\rm I}_{3}{\rm A}(F_{3}) \setminus {\rm I}_{4}{\rm A}(F_{3})$. Furthermore, it is easily shown that, for $a_{1} + a_{2} + a_{3} \neq 0$, $h \in \gamma_{2}(H)$ and $h \notin {\rm I}_{4}{\rm A}(F_{3})$. Therefore, $\gamma_{2}(H) = H \cap {\rm I}_{3}{\rm A}(F_{3})$. 

For $c = 3$, we apply similar arguments as before. We point out that $\gamma_{3}(H)/\gamma_{4}(H)$ is a free abelian group of rank $8$. Each element $g \in \gamma_{3}(H) \setminus \gamma_{4}(H)$ is uniquely written as $g = g_{1}g_{2}g_{3} u$, where $g_{1} = (z,x,x)^{b_{1}}(z,x,z)^{b_{2}}$, $g_{2} = (y,x,x)^{b_{3}} (y,x,z)^{b_{4}}(z,x,y)^{b_{5}} (z,y,z)^{b_{6}}$, $g_{3} = (y,x,y)^{b_{7}}(z,y,y)^{b_{8}}$, $u \in \gamma_{4}(H)$ and $b_{1}, \ldots, b_{8}$ are non-negative integers. Since $g_{1}(x_{i}) = x_{i}$ for $i \neq 1, 3$, we have, by direct calculations, 
$$
g_{1}(x_{2}) = x_{2}(x_{3},x_{1},x_{1},x_{2})^{-b_{1}}(x_{3},x_{1},x_{2},x_{3})^{-b_{2}} (x_{3},x_{2},(x_{3},x_{1}))^{b_{2}}v_{12}, \eqno(1)
$$
where $v_{12} \in \gamma_{5}(F_{3})$. By the equation (*), and working modulo $\gamma_{4}(H)$, we have 
$$
\begin{array}{ccc}
\tau^{b_{3}}_{(x_{2},x_{1},x_{1})} \equiv (y,x,x)^{b_{3}} (\chi_{32},x,x)^{b_{3}} & ~{\rm and}~& \tau^{b_{4}}_{(x_{3},x_{2},x_{1})} \equiv (y,x,z)^{b_{4}} (\chi_{32},x,z)^{b_{4}}.
\end{array}
$$
Since both $(\chi_{32},x,x)^{b_{3}}$ and $(\chi_{32},x,z)^{b_{4}}$ fix $x_{1}$, we get 
$$
\begin{array}{ccc}
(y,x,x)^{b_{3}}(x_{1}) = x_{1}(x_{1}, (x_{2}, x_{1}, x_{1}))^{b_{3}}v_{3} & ~{\rm and}~& (y,x,z)^{b_{4}}(x_{1}) = x_{1}(x_{1},(x_{3},x_{2},x_{1}))^{b_{4}}v_{4},
\end{array}
$$
where $v_{3}, v_{4} \in \gamma_{5}(F_{3})$. Since $(z,x)(x_{2}) = x_{2}(x_{3},x_{1},x_{2})^{-1}u$, with $u \in \gamma_{4}(F_{3})$, and by the equation (*), we get 
$$
(z,x,\tau_{x^{-1}_{2}})^{b_{5}} \equiv \tau^{b_{5}}_{(x_{3},x_{1},x_{2})^{-1}} \equiv (z,x,y)^{b_{5}} (z,x,\chi_{32})^{b_{5}}$$
and
$$(\tau_{(x_{3}^{-1},x_{2}^{-1})},\chi_{23})^{b_{6}} \equiv \tau^{b_{6}}_{(x_{3},x_{2},x^{-1}_{3})} \equiv (z,y,z)^{b_{6}} (z,y,\chi_{32})^{b_{6}}. 
$$
Since both $(z,x,\chi_{32})^{b_{5}}$ and $(z,y,\chi_{32})^{b_{6}}$ fix $x_{1}$, we obtain 
$$
\begin{array}{ccc}
(z,x,y)^{b_{5}}(x_{1}) = x_{1}(x_{3}, x_{1}, x_{2}, x_{1})^{b_{5}}v_{5} & ~{\rm and}~& (z,y,z)^{b_{6}}(x_{1}) = x_{1}(x_{3},x_{2},x_{3},x_{1})^{b_{6}}v_{6},
\end{array}
$$
where $v_{5}, v_{6} \in \gamma_{5}(F_{3})$. Therefore, 
$$
g_{2}(x_{1}) = x_{1}(x_{2}, x_{1}, x_{1},x_{1})^{-b_{3}}(x_{3},x_{2},x_{1},x_{1})^{-b_{4}}(x_{3},x_{1},x_{2},x_{1})^{b_{5}}(x_{3},x_{2},x_{3},x_{1})^{b_{6}}v, \eqno(2)
$$
where $v \in \gamma_{5}(F_{3})$. 

Next, we point out that $(\chi_{32},\chi_{21})(x_{i}) = x_{i}$ for $i = 1,2$ and $(\chi_{32},\chi_{21})(x_{3}) = x_{3}(x_{2},x_{1},x_{3})w$, where $w \in \gamma_{4}(F_{3})$. Since $(\chi_{32},\chi_{21}, \tau^{-1}_{x_{2}}) = {\rm Id}_{F_{3}}$, we have 
$$
(\chi_{32},x,y) = (\chi_{32},x,\chi_{32})^{-1} \psi,
$$
where $\psi \in \gamma_{4}(M_{3})$. Hence, 
$$
(\tau_{(\tau^{-1}_{2},\tau^{-1}_{1})},y)^{b_{7}} \equiv (y,x,y)^{b_{7}}(\chi_{32},x,\chi_{32})^{-b_{7}}
$$
and so,
$$
(\tau_{(x^{-1}_{2},x^{-1}_{1})},\tau_{x^{-1}_{2}}\chi^{-1}_{32})^{b_{7}} \equiv (y,x,y)^{b_{7}}(\chi_{32},x,\chi_{32})^{-b_{7}}.
$$
Since $(\tau_{(x^{-1}_{2},x^{-1}_{1})},\chi^{-1}_{32})^{b_{7}} = {\rm Id}_{F_{3}}$, we get
$$
(\tau_{(x^{-1}_{2},x^{-1}_{1})},\tau_{x^{-1}_{2}})^{b_{7}} \equiv (y,x,y)^{b_{7}}(\chi_{32},x,\chi_{32})^{-b_{7}}.
$$ 
Therefore,
$$
(y,x,y)^{b_{7}}(x_{1}) = x_{1}(x_{2},x_{1},x_{1},x_{2})^{b_{7}}v_{7},
$$
where $v_{7} \in \gamma_{5}(F_{3})$. By direct calculations, $(z,\chi_{32})(x_{1}) = x_{1}$, $(z,\chi_{32})(x_{2}) = x_{2}(x_{3},x_{2},x_{2})^{-1}u_{2}$ and $(z,\chi_{32})(x_{3}) = x_{3}(x_{3},x_{2},x_{3})u_{3}$, where $u_{2}, u_{3} \in \gamma_{4}(F_{3})$. It is easily verified that 
$$
(z,\chi_{32},\tau_{x^{-1}_{2}}) = \tau_{(x_{3},x_{2},x_{2})^{-1}} \psi_{1},
$$
where $\psi_{1} \in \gamma_{4}(M_{3})$. By the equation (*), 
$$
(\tau_{(x^{-1}_{3},x^{-1}_{2})},y)^{b_{8}} \equiv (z,y,y)^{b_{8}}(z,\chi_{32},y)^{b_{8}}.
$$
Since 
$$
(\tau_{(x^{-1}_{3},x^{-1}_{2})},y)^{b_{8}} \equiv (\tau_{(x^{-1}_{3},x^{-1}_{2})},\tau_{x^{-1}_{2}})^{b_{8}}(\tau_{(x^{-1}_{3},x^{-1}_{2})},\chi^{-1}_{32})^{b_{8}}
$$
and 
$$
(z,\chi_{32},y)^{b_{8}} \equiv (z,\chi_{32},\tau_{x^{-1}_{2}})^{b_{8}} (z,\chi_{32},\chi_{32})^{-b_{8}},
$$
we have 
$$
(\tau_{(x_{3},x_{2})},\chi^{-1}_{32})^{b_{8}} \equiv (z,y,y)^{b_{8}}(z,\chi_{32},\chi_{32})^{-b_{8}}.
$$
Hence, 
$$
\tau^{b_{8}}_{(x_{3},x_{2},x_{2})^{-1}} \equiv (z,y,y)^{b_{8}}(z,\chi_{32},\chi_{32})^{-b_{8}}.
$$
Since $(z,\chi_{32},\chi_{32})$ fixes $x_{1}$, we have 
$$
(z,y,y)^{b_{8}}(x_{1}) = x_{1}(x_{3},x_{2},x_{2},x_{1})^{b_{8}}v_{8},
$$
where $v_{8} \in \gamma_{5}(F_{3})$. Therefore,
$$
g_{3}(x_{1}) = x_{1}(x_{2},x_{1},x_{1},x_{2})^{b_{7}} (x_{3},x_{2},x_{2},x_{1})^{b_{8}}v_{7,8}, \eqno(3)
$$
where $v_{7,8} \in \gamma_{5}(F_{3})$. Since the "basic" group commutators of length $4$ consists of a basis for $\gamma_{4}(F_{3})/\gamma_{5}(F_{3})$, we obtain from the equations (1), (2) and (3) that $g \in \gamma_{4}(H)$ and $g \notin {\rm I}_{6}{\rm A}(F_{3})$, that is, $\gamma_{4}(H) = H \cap {\rm I}_{5}{\rm A}(F_{3})$.

}
\end{remark}

\bigskip

\noindent V. Metaftsis, Department of Mathematics, University of
the Aegean, Karlovassi, 832 00 Samos, Greece. {\it e-mail:}
vmet@aegean.gr

\bigskip

\noindent A.I. Papistas, Department of Mathematics, Aristotle
University of Thessaloniki, 541 24 Thessaloniki, Greece. {\it
e-mail:} apapist@math.auth.gr

\bigskip

\noindent H. Sevaslidou, Department of Mathematics, Aristotle
University of Thessaloniki, 541 24 Thessaloniki, Greece. {\it
e-mail:} sashasevaslidou@gmail.com

\end{document}